\documentstyle[12pt]{article}
\newcommand{\binom}[2]{\mbox{$\left( #1 \atop #2 \right)$}}

\newtheorem{thrm}{Theorem}

\begin{document}
\title{GRAPH COMPOSITIONS I:\\Basic enumeration}
\author{A. Knopfmacher\\{\it University of the Witwatersrand}
\and
M. E. Mays\\{\it West Virginia University}
\thanks 
{The second author thanks the Centre for 
Applicable Analysis and Number Theory at The University of
the Witwatersrand for sponsoring his visit during May and
June 1998.}}
\date{February 3, 2000}
\maketitle
\begin{abstract} 
The idea of graph compositions generalizes both ordinary compositions 
of positive integers and partitions of finite sets.  In this paper we 
develop formulas, generating functions, and recurrence relations for 
composition counting functions for several families of graphs.
\end{abstract}

\section{Introduction}

Let $G$ be a labelled graph, with edge set $E(G)$ and vertex set 
$V(G)$.  A {\em composition} of $G$ is a partition of $V(G)$ into
vertex sets of connected induced subgraphs of $G$. 
Thus a partition provides a set of connected subgraphs of $G$, 
$\{G_1,G_2,\cdots,G_m\}$, 
with the properties that 
$\bigcup_{i=1}^m V(G_i) = V(G)$ and
for $i \neq j, V(G_i) \cap V(G_j) = \emptyset$.  (Note, however,
that since different edge subsets of a graph can span the same
vertex set, it is possible for a different set of connected subgraphs
of $G$ to yield the same composition.)
We will call the vertex sets $V(G_i)$, or the subgraphs $G_i$ 
themselves if there is no danger of confusion, {\em components} 
of a given 
composition.  This paper is most concerned with straightforward
enumerative questions:
counting how many compositions a given graph has.  Topics such
as restricted compositions or asymptotic results will be considered
later.  We will denote by
$C(G)$ the number of distinct compositions that exist for a given
graph $G$.

For example, the complete bipartite graph $K_{2,3}$ has exactly
34 compositions, which are illustrated below.  The significance
of the edges shown is to indicate the connected components:  it is possible
that other choices of edges could yield the same connected components,
and hence the same composition.  In fact,
since there are 64 subsets of the set of six edges of $K_{2,3}$, this 
overlap must occur.

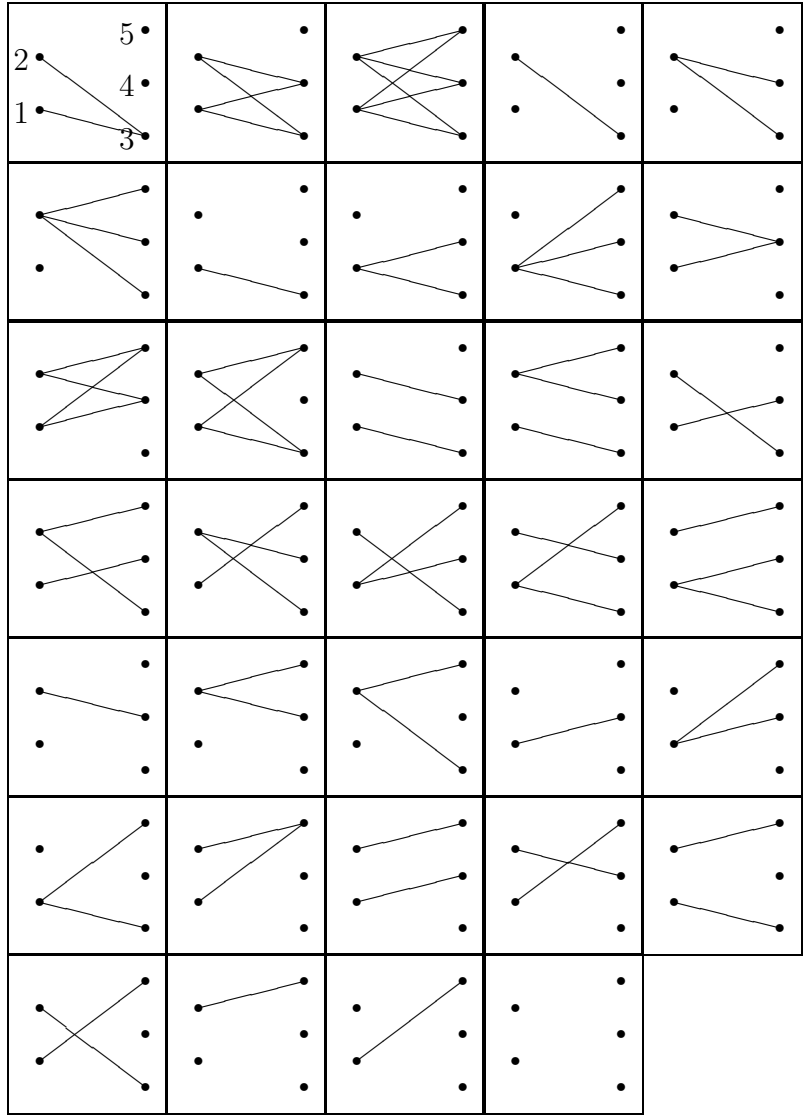
\begin{picture}(300,420)
%1
\put(0,360){\framebox(60,60){
\begin{picture}(60, 60)
\put(0,15){1}
\put(0,35){2}
\put(40,5){3}
\put(40,25){4}
\put(40,45){5}
\put(10,20){\circle*{3}}
\put(10,40){\circle*{3}}
\put(50,10){\circle*{3}}
\put(50,30){\circle*{3}}
\put(50,50){\circle*{3}}
\put(10,20){\line(4,-1){40}}
\put(10,40){\line(4,-3){40}}
\end{picture}}}
%2
\put(60,360){\framebox(60,60){
\begin{picture}(60, 60)
\put(10,20){\circle*{3}}
\put(10,40){\circle*{3}}
\put(50,10){\circle*{3}}
\put(50,30){\circle*{3}}
\put(50,50){\circle*{3}}
\put(10,20){\line(4,1){40}}
\put(10,20){\line(4,-1){40}}
\put(10,40){\line(4,-1){40}}
\put(10,40){\line(4,-3){40}}
\end{picture}}}
%3
\put(120,360){\framebox(60,60){
\begin{picture}(60, 60)
\put(10,20){\circle*{3}}
\put(10,40){\circle*{3}}
\put(50,10){\circle*{3}}
\put(50,30){\circle*{3}}
\put(50,50){\circle*{3}}
\put(10,20){\line(4,3){40}}
\put(10,20){\line(4,1){40}}
\put(10,20){\line(4,-1){40}}
\put(10,40){\line(4,1){40}}
\put(10,40){\line(4,-1){40}}
\put(10,40){\line(4,-3){40}}
\end{picture}}}
%4
\put(180,360){\framebox(60,60){
\begin{picture}(60, 60)
\put(10,20){\circle*{3}}
\put(10,40){\circle*{3}}
\put(50,10){\circle*{3}}
\put(50,30){\circle*{3}}
\put(50,50){\circle*{3}}
\put(10,40){\line(4,-3){40}}
\end{picture}}}
%5
\put(240,360){\framebox(60,60){
\begin{picture}(60, 60)
\put(10,20){\circle*{3}}
\put(10,40){\circle*{3}}
\put(50,10){\circle*{3}}
\put(50,30){\circle*{3}}
\put(50,50){\circle*{3}}
\put(10,40){\line(4,-1){40}}
\put(10,40){\line(4,-3){40}}
\end{picture}}}
%6
\put(0,300){\framebox(60,60){
\begin{picture}(60, 60)
\put(10,20){\circle*{3}}
\put(10,40){\circle*{3}}
\put(50,10){\circle*{3}}
\put(50,30){\circle*{3}}
\put(50,50){\circle*{3}}
\put(10,40){\line(4,1){40}}
\put(10,40){\line(4,-1){40}}
\put(10,40){\line(4,-3){40}}
\end{picture}}}
%7
\put(60,300){\framebox(60,60){
\begin{picture}(60, 60)
\put(10,20){\circle*{3}}
\put(10,40){\circle*{3}}
\put(50,10){\circle*{3}}
\put(50,30){\circle*{3}}
\put(50,50){\circle*{3}}
\put(10,20){\line(4,-1){40}}
\end{picture}}}
%8
\put(120,300){\framebox(60,60){
\begin{picture}(60, 60)
\put(10,20){\circle*{3}}
\put(10,40){\circle*{3}}
\put(50,10){\circle*{3}}
\put(50,30){\circle*{3}}
\put(50,50){\circle*{3}}
\put(10,20){\line(4,1){40}}
\put(10,20){\line(4,-1){40}}
\end{picture}}}
%9
\put(180,300){\framebox(60,60){
\begin{picture}(60, 60)
\put(10,20){\circle*{3}}
\put(10,40){\circle*{3}}
\put(50,10){\circle*{3}}
\put(50,30){\circle*{3}}
\put(50,50){\circle*{3}}
\put(10,20){\line(4,3){40}}
\put(10,20){\line(4,1){40}}
\put(10,20){\line(4,-1){40}}
\end{picture}}}
%10
\put(240,300){\framebox(60,60){
\begin{picture}(60, 60)
\put(10,20){\circle*{3}}
\put(10,40){\circle*{3}}
\put(50,10){\circle*{3}}
\put(50,30){\circle*{3}}
\put(50,50){\circle*{3}}
\put(10,20){\line(4,1){40}}
\put(10,40){\line(4,-1){40}}
\end{picture}}}
%11
\put(0,240){\framebox(60,60){
\begin{picture}(60, 60)
\put(10,20){\circle*{3}}
\put(10,40){\circle*{3}}
\put(50,10){\circle*{3}}
\put(50,30){\circle*{3}}
\put(50,50){\circle*{3}}
\put(10,20){\line(4,3){40}}
\put(10,20){\line(4,1){40}}
\put(10,40){\line(4,1){40}}
\put(10,40){\line(4,-1){40}}
\end{picture}}}
%12
\put(60,240){\framebox(60,60){
\begin{picture}(60, 60)
\put(10,20){\circle*{3}}
\put(10,40){\circle*{3}}
\put(50,10){\circle*{3}}
\put(50,30){\circle*{3}}
\put(50,50){\circle*{3}}
\put(10,20){\line(4,3){40}}
\put(10,20){\line(4,-1){40}}
\put(10,40){\line(4,1){40}}
\put(10,40){\line(4,-3){40}}
\end{picture}}}
%13
\put(120,240){\framebox(60,60){
\begin{picture}(60, 60)
\put(10,20){\circle*{3}}
\put(10,40){\circle*{3}}
\put(50,10){\circle*{3}}
\put(50,30){\circle*{3}}
\put(50,50){\circle*{3}}
\put(10,20){\line(4,-1){40}}
\put(10,40){\line(4,-1){40}}
\end{picture}}}
%14
\put(180,240){\framebox(60,60){
\begin{picture}(60, 60)
\put(10,20){\circle*{3}}
\put(10,40){\circle*{3}}
\put(50,10){\circle*{3}}
\put(50,30){\circle*{3}}
\put(50,50){\circle*{3}}
\put(10,20){\line(4,-1){40}}
\put(10,40){\line(4,1){40}}
\put(10,40){\line(4,-1){40}}
\end{picture}}}
%15
\put(240,240){\framebox(60,60){
\begin{picture}(60, 60)
\put(10,20){\circle*{3}}
\put(10,40){\circle*{3}}
\put(50,10){\circle*{3}}
\put(50,30){\circle*{3}}
\put(50,50){\circle*{3}}
\put(10,20){\line(4,1){40}}
\put(10,40){\line(4,-3){40}}
\end{picture}}}
%16
\put(0,180){\framebox(60,60){
\begin{picture}(60, 60)
\put(10,20){\circle*{3}}
\put(10,40){\circle*{3}}
\put(50,10){\circle*{3}}
\put(50,30){\circle*{3}}
\put(50,50){\circle*{3}}
\put(10,20){\line(4,1){40}}
\put(10,40){\line(4,1){40}}
\put(10,40){\line(4,-3){40}}
\end{picture}}}
%17
\put(60,180){\framebox(60,60){
\begin{picture}(60, 60)
\put(10,20){\circle*{3}}
\put(10,40){\circle*{3}}
\put(50,10){\circle*{3}}
\put(50,30){\circle*{3}}
\put(50,50){\circle*{3}}
\put(10,20){\line(4,3){40}}
\put(10,40){\line(4,-1){40}}
\put(10,40){\line(4,-3){40}}
\end{picture}}}
%18
\put(120,180){\framebox(60,60){
\begin{picture}(60, 60)
\put(10,20){\circle*{3}}
\put(10,40){\circle*{3}}
\put(50,10){\circle*{3}}
\put(50,30){\circle*{3}}
\put(50,50){\circle*{3}}
\put(10,20){\line(4,3){40}}
\put(10,20){\line(4,1){40}}
\put(10,40){\line(4,-3){40}}
\end{picture}}}
%19
\put(180,180){\framebox(60,60){
\begin{picture}(60, 60)
\put(10,20){\circle*{3}}
\put(10,40){\circle*{3}}
\put(50,10){\circle*{3}}
\put(50,30){\circle*{3}}
\put(50,50){\circle*{3}}
\put(10,20){\line(4,3){40}}
\put(10,20){\line(4,-1){40}}
\put(10,40){\line(4,-1){40}}
\end{picture}}}
%20
\put(240,180){\framebox(60,60){
\begin{picture}(60, 60)
\put(10,20){\circle*{3}}
\put(10,40){\circle*{3}}
\put(50,10){\circle*{3}}
\put(50,30){\circle*{3}}
\put(50,50){\circle*{3}}
\put(10,20){\line(4,1){40}}
\put(10,20){\line(4,-1){40}}
\put(10,40){\line(4,1){40}}
\end{picture}}}
%21
\put(0,120){\framebox(60,60){
\begin{picture}(60, 60)
\put(10,20){\circle*{3}}
\put(10,40){\circle*{3}}
\put(50,10){\circle*{3}}
\put(50,30){\circle*{3}}
\put(50,50){\circle*{3}}
\put(10,40){\line(4,-1){40}}
\end{picture}}}
%22
\put(60,120){\framebox(60,60){
\begin{picture}(60, 60)
\put(10,20){\circle*{3}}
\put(10,40){\circle*{3}}
\put(50,10){\circle*{3}}
\put(50,30){\circle*{3}}
\put(50,50){\circle*{3}}
\put(10,40){\line(4,1){40}}
\put(10,40){\line(4,-1){40}}
\end{picture}}}
%23
\put(120,120){\framebox(60,60){
\begin{picture}(60, 60)
\put(10,20){\circle*{3}}
\put(10,40){\circle*{3}}
\put(50,10){\circle*{3}}
\put(50,30){\circle*{3}}
\put(50,50){\circle*{3}}
\put(10,40){\line(4,1){40}}
\put(10,40){\line(4,-3){40}}
\end{picture}}}
%24
\put(180,120){\framebox(60,60){
\begin{picture}(60, 60)
\put(10,20){\circle*{3}}
\put(10,40){\circle*{3}}
\put(50,10){\circle*{3}}
\put(50,30){\circle*{3}}
\put(50,50){\circle*{3}}
\put(10,20){\line(4,1){40}}
\end{picture}}}
%25
\put(240,120){\framebox(60,60){
\begin{picture}(60, 60)
\put(10,20){\circle*{3}}
\put(10,40){\circle*{3}}
\put(50,10){\circle*{3}}
\put(50,30){\circle*{3}}
\put(50,50){\circle*{3}}
\put(10,20){\line(4,3){40}}
\put(10,20){\line(4,1){40}}
\end{picture}}}
%26
\put(0,60){\framebox(60,60){
\begin{picture}(60, 60)
\put(10,20){\circle*{3}}
\put(10,40){\circle*{3}}
\put(50,10){\circle*{3}}
\put(50,30){\circle*{3}}
\put(50,50){\circle*{3}}
\put(10,20){\line(4,3){40}}
\put(10,20){\line(4,-1){40}}
\end{picture}}}
%27
\put(60,60){\framebox(60,60){
\begin{picture}(60, 60)
\put(10,20){\circle*{3}}
\put(10,40){\circle*{3}}
\put(50,10){\circle*{3}}
\put(50,30){\circle*{3}}
\put(50,50){\circle*{3}}
\put(10,20){\line(4,3){40}}
\put(10,40){\line(4,1){40}}
\end{picture}}}
%28
\put(120,60){\framebox(60,60){
\begin{picture}(60, 60)
\put(10,20){\circle*{3}}
\put(10,40){\circle*{3}}
\put(50,10){\circle*{3}}
\put(50,30){\circle*{3}}
\put(50,50){\circle*{3}}
\put(10,20){\line(4,1){40}}
\put(10,40){\line(4,1){40}}
\end{picture}}}
%29
\put(180,60){\framebox(60,60){
\begin{picture}(60, 60)
\put(10,20){\circle*{3}}
\put(10,40){\circle*{3}}
\put(50,10){\circle*{3}}
\put(50,30){\circle*{3}}
\put(50,50){\circle*{3}}
\put(10,20){\line(4,3){40}}
\put(10,40){\line(4,-1){40}}
\end{picture}}}
%30
\put(240,60){\framebox(60,60){
\begin{picture}(60, 60)
\put(10,20){\circle*{3}}
\put(10,40){\circle*{3}}
\put(50,10){\circle*{3}}
\put(50,30){\circle*{3}}
\put(50,50){\circle*{3}}
\put(10,20){\line(4,-1){40}}
\put(10,40){\line(4,1){40}}
\end{picture}}}
%31
\put(0,0){\framebox(60,60){
\begin{picture}(60, 60)
\put(10,20){\circle*{3}}
\put(10,40){\circle*{3}}
\put(50,10){\circle*{3}}
\put(50,30){\circle*{3}}
\put(50,50){\circle*{3}}
\put(10,20){\line(4,3){40}}
\put(10,40){\line(4,-3){40}}
\end{picture}}}
%32
\put(60,0){\framebox(60,60){
\begin{picture}(60, 60)
\put(10,20){\circle*{3}}
\put(10,40){\circle*{3}}
\put(50,10){\circle*{3}}
\put(50,30){\circle*{3}}
\put(50,50){\circle*{3}}
\put(10,40){\line(4,1){40}}
\end{picture}}}
%33
\put(120,0){\framebox(60,60){
\begin{picture}(60, 60)
\put(10,20){\circle*{3}}
\put(10,40){\circle*{3}}
\put(50,10){\circle*{3}}
\put(50,30){\circle*{3}}
\put(50,50){\circle*{3}}
\put(10,20){\line(4,3){40}}
\end{picture}}}
%34
\put(180,0){\framebox(60,60){
\begin{picture}(60, 60)
\put(10,20){\circle*{3}}
\put(10,40){\circle*{3}}
\put(50,10){\circle*{3}}
\put(50,30){\circle*{3}}
\put(50,50){\circle*{3}}
\end{picture}}}
\end{picture}

Theorem~\ref{path} below is a well known result that motivates this choice of 
terminology, and Theorem~\ref{complete} relates 
the idea to another familiar combinatorial setting.

Let $G=P_n$, the path with $n$ vertices.  Then any 
subgraph of $G$ is also a path, and the components of a composition
consist of paths of cardinality $|G_i|=a_i$ so that 
$\sum_{i=1}^m a_i = n$.  Thus the path lengths provide a composition
of the positive integer $n$ (a representation of $n$ as an ordered sum
of positive integers), and any composition of $n$ determines ``cut
points'' to provide a composition of the graph $P_n$.  The well known   
counting function for integer compositions applies to give the first
result.

\begin{thrm}\label{path}
$C(P_n)=2^{n-1}$.
\end{thrm}

We will define $C(P_0)$ to be 1 in order to make a formula in
Theorem~\ref{wheel} below more palatable.

Now we consider another case, a family of graphs with many edges.  
Let $G=K_n$, the complete 
graph on $n$ vertices.  Then any
subset of $V(G)$ can serve as the vertex set of a subgraph of $G$, and
the number of compositions of $G$ is the number of partitions of a set 
with $n$ elements into nonempty subsets.  The number of partitions of a set
of $n$ elements is given by the Bell number $B(n)$.  The sequence of Bell
numbers begins $1, 2, 5 ,15, 52, \cdots$, and has generating function
$e^{e^x-1}$.  
This well known sequence has an
extensive bibliography compiled by Gould \cite{gould}.  

\begin{thrm}\label{complete} 
$C(K_n)=B(n)$.
\end{thrm}

These two results are extreme cases:  no connected
graph $G$ with $n$ vertices can have fewer than $C(P_n)$ 
compositions, nor more than $C(K_n)$.  Thus for $\{F_n\}_{n \ge 1}$
a family of connected graphs such that $|V(F_n)|=n$, the
values $C(F_n)$ satisfy $2^{n-1} \le C(F_n) \le B(n)$.
We allow graphs to be disconnected, and the extreme case would
be the graph with no edges, and $n$ isolated vertices.  By our
definition this graph has exactly one composition.
 
\section{General observations}
  In general, one might expect that for graphs with a given number of
vertices, the more edges, the more compositions.  This is not always true,
and certainly more information is needed than $|V(G)|$ and
$|E(G)|$
to determine $C(G)$.  The example below shows two graphs $G_1$ 
and $G_2$ with 4 vertices and 4 edges, but 
$C(G_1)=10 \neq 12=C(G_2)$.

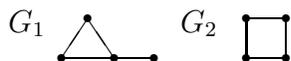
\begin{picture}(260,60)
\put(30,10){\circle*{3}}
\put(50,10){\circle*{3}}
\put(40,25){\circle*{3}}
\put(65,10){\circle*{3}}
\put(30,10){\line(1,0){20}}
\put(30,10){\line(2,3){10}}
\put(50,10){\line(-2,3){10}}
\put(50,10){\line(1,0){15}}
\put(10,20){$G_1$}
\put(100,10){\circle*{3}}
\put(115,10){\circle*{3}}
\put(100,25){\circle*{3}}
\put(115,25){\circle*{3}}
\put(100,10){\line(1,0){15}}
\put(100,10){\line(0,1){15}}
\put(100,25){\line(1,0){15}}
\put(115,10){\line(0,1){15}}
\put(75,20){$G_2$}
\end{picture}

\begin{thrm}\label{disconnected}
If $G=G_1 \cup G_2$ and there are no edges from vertices in $G_1$
to vertices in $G_2$ (i.e. $G$ is {\em disconnected}), then 
$C(G)=C(G_1) \cdot C(G_2)$.  The same result holds if $G_1$ and $G_2$
have exactly one vertex in common.
\end{thrm}
{\sc Proof.} This is a consequence of the Fundamental Principle 
of Counting.  We obtain compositions of $G$ by pairing compositions
of $G_1$ with compositions of $G_2$ in all possible ways. $\Box$

We can also give a general result for graphs that are 
``almost disconnected''.

\begin{thrm}\label{one edge}
If $G=G_1 \cup G_2$ and there is an edge from one of the 
vertices of $G_1$ to one of the vertices of $G_2$ whose
removal disconnects $G$, then 
$C(G)=2 \cdot C(G_1) \cdot C(G_2)$.
\end{thrm}
{\sc Proof.} Call the distinguished edge $e$, between vertices $v_i$
and $v_j$. For any composition
of $G_1$ and any composition of $G_2$ we can build a composition
of $G$ in exactly two ways:  either $e$ can be included to combine
the component of $v_i$ in $G_1$ and the component of $v_j$ in $G_2$,
or not.  Thus the count provided by Theorem~\ref{disconnected}
is doubled. $\Box$

The analysis when $G$ consists of two subgraphs connected by a
bridge of
$n>1$ vertices is more complicated.  More information is required
about the nature of the components containing the connecting
vertices in compositions of the subgraphs.  Several special cases
are considered in later sections.

\begin{thrm}\label{tree}
Let $T_n$ be any tree with $n$ vertices.  Then $C(T_n)=2^{n-1}$.
\end{thrm}
{\sc Proof.} The proof is by induction.  When $n=1$ the tree is a
single vertex, with $1=2^0$ compositions.  If the result is true
for $n \le k$, we consider $T_{k+1}$ and remove an edge.  This 
disconnects $T_{k+1}$, into two subtrees with $l$ and $k+1-l$
vertices for some $l \ge 1$.  The induction hypothesis applies
to each subtree, giving $2^{l-1}$ and $2^{k-l}$ compositions.
Theorem~\ref{one edge} then gives
$2 \cdot 2^{l-1} \cdot 2^{k-l}=2^{k}$
compositions for $T_{k+1}$. $\Box$

The {\em star} graph $S_n$ consists of a distinguished {\em center}
vertex connected to each of $n-1$ {\em edge} vertices.  $S_n$ is an
example of a tree, and so $C(S_n)=2^{n-1}$.

Deleting one edge from a complete graph has a predictable effect.

\begin{thrm}\label{complete-1}
Let $K_n^-$ denote the complete graph on $n$ vertices with one edge
removed.  Then C($K_n^-) = B(n)-B(n-2)$.
\end{thrm}
{\sc Proof.} The only time that the deleted edge $e$ between $v_i$
and $v_j$ affects a composition
counted by $C(K_n)$ is when the component containing $v_i$ and $v_j$
consists of exactly those two vertices.  Otherwise there is a path
between $v_i$ and $v_j$ in $K_n$ bypassing the deleted edge.  Hence
from the $B(n)$ compositions counted by $C(K_n)$ must be deleted
exactly those compositions for which one component is $\{v_i,v_j\}$.
This restriction rules out exactly $C(K_{n-2})=B(n-2)$ compositions
of $K_n$. $\Box$

On the other hand, deleting more than one edge affects the number of
compositions
depending on whether the edges deleted are adjacent or not.  For example
the graph resulting when two adjacent edges are deleted from $K_5$
has 40 compositions, whereas if two nonadjacent edges are deleted
the resulting graph has 43 compositions.

Another basic family of graphs to consider are the cycle graphs $C_n$.
$C_n$ is the graph with $n$ vertices and $n$ edges, with vertex $i$
connected to vertices $i \pm 1 \pmod{n}$.

\begin{thrm}\label{cycle} 
$C(C_n)=2^n-n$
\end{thrm}
{\sc Proof.} Pick any edge of the cycle and delete it.  The resulting
graph is $P_n$, with $C(P_n)=2^{n-1}$ by Theorem~\ref{path}.  Any composition of $P_n$ may
be regarded as a composition of $C_n$ as well.  The deleted edge
may be reinserted, providing a new composition of $C_n$ not previously
counted, unless the composition of $P_n$ had been obtained by deleting
no edge, or exactly one edge, from $P_n$.  In these cases,
reinserting the original deleted edge results in the same composition
of $C_n$: the composition consisting of the single component consisting
of all $n$ vertices.  Hence the total count of distinct compositions
of $C_n$ is $2 \cdot 2^{n-1}-n=2^n-n$. $\Box$  

It is sometimes useful to group the compositions of $C_n$ so that 
different compositions obtained by rotation may be analysed together.
This idea has its origins in the general area of combinatorics on
words, where periodicity and cyclic permutations are studied via
what are called Lyndon words \cite{cm}, \cite{mays}.
Analogously, we define a {\em Lyndon composition} of the positive
integer $n$ to be an aperiodic composition that is lexicographically
least among its cyclic permutations.  For example, $1+2+1+2$ is not
a Lyndon composition of 6 because it is periodic, and $1+1+2+2$ is
a Lyndon composition of 6 because it is aperiodic, and in addition 
by the lexicographic ordering we order the cyclic permutations of
the summands as
``1+1+2+2'' $<$ ``1+2+2+1'' $<$ ``2+1+1+2'' $<$ ``2+2+1+1''.
The number of Lyndon compositions $L(n)$ of the integer $n$ is given 
by the formula
\begin{equation}
L(n)=\frac{1}{n} \sum_{d \mid n} \mu(\frac{n}{d}) 2^d. \label{lyndon} 
\end{equation}
By (\ref{lyndon}) we should define $L(1)=2$.  Then 
\begin{displaymath}
C(C_n)=\sum_{d \mid n} d L(d) - n,
\end{displaymath}
which, together with the inverted version of (\ref{lyndon}), 
recovers the formula in Theorem~\ref{cycle}.  We will have use for the
sequence of values of $L(n)$:
\[2, 1, 2, 3, 6, 9, 18, 30, 56, \ldots \]
 
The wheel graph $W_n$ consists of the star graph $S_n$ with extra edges
appended so that there is a cycle through the $n-1$ outer vertices.
Alternately, $W_n$ is $C_{n-1}$
with one extra ``central'' vertex appended which
is adjacent to each ``outer'' vertex in the cycle.  We will take $W_1$ to be
an isolated single vertex, $W_2$ to be $P_2$, and $W_3$ to be
$C_3$.  Then the sequence
$\{C(W_n)\}$ begins
\[1, 2, 5, 15, 43, 118, 316, 836, 2199, 5769, 15117, 39592, \ldots.\]
We account for these values in the theorem below.

\begin{thrm}\label{wheel}
\begin{displaymath}
C(W_n)= 2^{n-1}-n+2+ \sum_{1<d \mid n-1} d
\sum_{a_1+\ldots+a_k=d}^{\hspace{.22in}\prime}
\prod_{i=1}^k  C(P_{a_i-1})^{(n-1)/d},
\end{displaymath}
where $\Sigma^\prime$ indicates a sum over Lyndon compositions of $d$.
\end{thrm}
{\sc Proof.}  There are two cases to consider.  Suppose first that
in a composition of $W_n$
the central vertex is connected to no outer vertex.  Then the
outer vertices may be grouped into $C(C_{n-1})=2^{n-1}-(n-1)$ distinct
compositions.  Now suppose that the central vertex is connected to
one or more outer vertices.  Then the remaining outer vertices are
disconnected into a set of paths.  The possible patterns of paths
correspond to Lyndon compositions of $n-1$ if they are not periodic,
or to adjoined Lyndon compositions of $d | n-1$ if they are periodic.
The correspondence is determined by using the number of gaps between
adjacent spokes of the wheel to be summands of the composition. The
number of compositions in this case is the product of the number
of compositions of the constituent paths.  This is the product term
in the summation formula.  The exponent of $(n-1)/d$ allows for
all possible combinations of paths in the case where there are
adjoined Lyndon compositions of proper divisors $d | n-1$.
$\Box$

We thank {\tt superseeker@research.att.com} for the observation
that the sequence of values of $C(W(n))$ corresponds to the
third difference of the bisection of the Lucas sequence. It 
also satisfies the recurrence
relation $C(W_1)=C(W_2)=2, C(W_n)=3C(W_{n-1})-C(W_{n-2})+n-2$.
There must be a combinatorial interpretation of this recurrence.

\section{Ladders $L_n$}
We build the ladder $L_n$ as a product of a path of length 2 and a
path of length $n$.  Thus $L_n$ has $2n$ vertices and $3n-2$ edges.
The four ``corner'' vertices have degree 2, and the other vertices have
degree 3.  We will take $L_1=P_2$, so $C(L_1)=2$.  $L_2=C_4$, so
$C(L_2)=12$ by Theorem \ref{cycle}.  The most direct way to account
for other values of $C(L_n)$ is with a recurrence.

\begin{thrm}\label{ladder} $C(L_1)=2, C(L_2)=12$, and for $n > 2,
C(L_n)=6 \cdot C(L_{n-1})+C(L_{n-2})$.
\end{thrm}
{\sc Proof.}  Label the vertices of $L_n$ as $a_{1,1}, a_{1,2},
a_{2,1}, a_{2,2}, \ldots, a_{n,1}, a_{n,2}$.  Denote by
$A_k$ the number of compositions of $L_k$ in which the vertices
$a_{n,1}$ and $a_{n,2}$ are in different components, and by $B_k$
the number of compositions of $L_k$ in which the vertices
$a_{n,1}$ and $a_{n,2}$ are in the same component.

In order to generate a composition of $L_n$ from $L_{n-1}$ there
are eight configurations to consider:

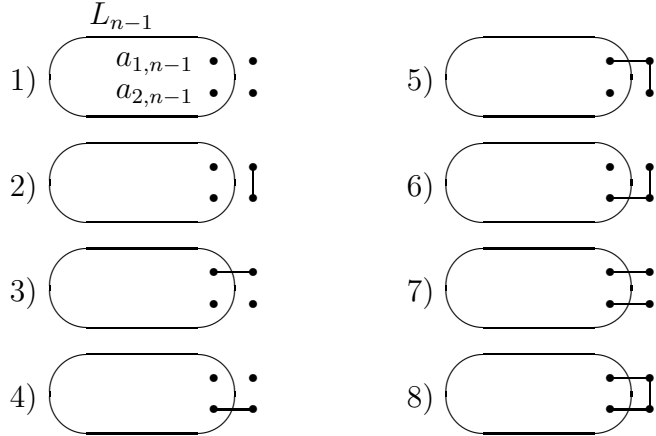
\begin{picture}(260,160)
\put(30,160){$L_{n-1}$}
\put(40,144){$a_{1,n-1}$}
\put(40,132){$a_{2,n-1}$}
\put(0,135){1)}
\put(77,146){\circle*{3}}
\put(77,134){\circle*{3}}
\put(92,146){\circle*{3}}
\put(92,134){\circle*{3}}
\put(50,140){\oval(70,30)}
\put(150,135){5)}
\put(227,146){\circle*{3}}
\put(227,134){\circle*{3}}
\put(242,146){\circle*{3}}
\put(242,134){\circle*{3}}
\put(200,140){\oval(70,30)}
\put(227,146){\line(1,0){15}}
\put(242,134){\line(0,1){12}}
\put(0,95){2)}
\put(77,106){\circle*{3}}
\put(77,94){\circle*{3}}
\put(92,106){\circle*{3}}
\put(92,94){\circle*{3}}
\put(50,100){\oval(70,30)}
\put(92,94){\line(0,1){12}}
\put(150,95){6)}
\put(227,106){\circle*{3}}
\put(227,94){\circle*{3}}
\put(242,106){\circle*{3}}
\put(242,94){\circle*{3}}
\put(200,100){\oval(70,30)}
\put(242,94){\line(0,1){12}}
\put(227,94){\line(1,0){15}}
\put(0,55){3)}
\put(77,66){\circle*{3}}
\put(77,54){\circle*{3}}
\put(92,66){\circle*{3}}
\put(92,54){\circle*{3}}
\put(50,60){\oval(70,30)}
\put(77,66){\line(1,0){15}}
\put(150,55){7)}
\put(227,66){\circle*{3}}
\put(227,54){\circle*{3}}
\put(242,66){\circle*{3}}
\put(242,54){\circle*{3}}
\put(200,60){\oval(70,30)}
\put(227,54){\line(1,0){15}}
\put(227,66){\line(1,0){15}}
\put(0,15){4)}
\put(77,26){\circle*{3}}
\put(77,14){\circle*{3}}
\put(92,26){\circle*{3}}
\put(92,14){\circle*{3}}
\put(50,20){\oval(70,30)}
\put(77,14){\line(1,0){15}}
\put(150,15){8)}
\put(227,26){\circle*{3}}
\put(227,14){\circle*{3}}
\put(242,26){\circle*{3}}
\put(242,14){\circle*{3}}
\put(200,20){\oval(70,30)}
\put(227,14){\line(1,0){15}}
\put(227,26){\line(1,0){15}}
\put(242,14){\line(0,1){12}}
\end{picture}

If we start with a composition counted by $A_{n-1}$, cases 1), 3),
4), and 7) yield distinct compositions counted by $A_n$.  If we start
with one counted by $B_{n-1}$, only 1), 3), and 4) yield distinct
compositions counted by $A_n$.  Hence
$A_n=4 \cdot A_{n-1} + 3 \cdot B_{n-1}$.
Similarly, cases 2), 5), and 6) go from a composition counted by
$A_{n-1}$ to one counted by $B_n$.  Starting with $B_{n-1}$, only
two distinct compositions arise:  the one given by case 2), or the
single new composition represented by cases 5), 6), 7) or 8).
Hence $B_n=3 \cdot A_{n-1} + 2 \cdot B_{n-1}$.  Since
$C(L_n)=A_n+B_n$, we have 
\[C(L_n)=7 \cdot A_{n-1} + 5 \cdot B_{n-1}.\]
On the other hand, 
\[A_{n-1}-B_{n-1}=A_{n-2}+B_{n-2}= C(L_{n-2}).\]  
Hence
\[C(L_n)=
6(A_{n-1}+B_{n-1})+(A_{n-1}-B_{n-1})=
6 \cdot C(L_{n-1})+C(L_{n-2}).\Box \]

As a bit of moonshine, we note that this recurrence guarantees the
sequence of values of $L_n/2$ matches the denominators in the continued
fraction expansion of $\sqrt{10}$.  A proof, but not an explanation,
is provided by observing recurrences and starting values are the same
for the two sequences.

\section{Bipartite graphs $K_{m,n}$}

An example showing that $C(K_{2,3})=34$ by exhibiting all 34
compositions
is in the first section.  The graphs
$K_{m,n}$, with $m+n$ vertices and $mn$ edges, are the most
complicated we will analyse in this paper.

\begin{thrm}\label{Kmn}
Define an array $A=(a_{i,j})$ via the recurrences
$a_{m,0}=0$ for any nonnegative integer $m$, $a_{0,1}=1,$
$a_{0,n}=0$ for any $n > 1$, and otherwise
\begin{equation}
a_{m,n}=\sum_{i=0}^{m-1} \binom{m-1}{i} a_{m-1-i,n-1} -
\sum_{i=1}^{m-1} \binom{m-1}{i} a_{m-1-i,n}. \label{a}
\end{equation}
Then
\begin{equation}
C(K_{m,n})= \sum_{i=1}^{m+1} a_{m,i} \, i^n. \label{Ksum}
\end{equation}
\end{thrm}
{\sc Proof.}
We observe $C(K_{m,0})=C(K_{0,n})=1$, vacuously.
$C(K_{m,1})=2^m$ because $K_{m,1}=S_{m+1}$, and similarly
for $K_{1,n}$.  This observation is the first step in an
induction on the arithmetic nature of 
$C(K_{m,n})$.  Now consider $C(K_{m,n})$ for $m \ge 1$.
Write the two parts of the bipartition as 
$A=\{a_1,a_2,\ldots,a_m\}$ and
$B=\{b_1,b_2,\ldots,b_n\}$.  $a_1$ must be in some component.
Consider cases.
\begin{trivlist}
\item[1)] $a_1$ is a singleton.  Then all the other components
determine a composition of $K_{m-1,n}$.  This can be done in
$C(K_{m-1,n})$ ways.
\item[2)] $a_1$ is in a component with no other elements of $A$,
but with elements of $B$.  Say a $j$-set of $B$.  The remaining
elements of $A$ and the remaining elements of $B$ can be paired
in $C(K_{m-1,n-j})$ ways.  There are $\binom{n}{j}$ $j$-sets of
$B$, so the total number of compositions here is 
\[
\sum_{j=1}^n \binom{n}{j} C(K_{m-1,n-j}).
\]
Cases 1) and 2) can be combined in a single sum:
\[
\sum_{j=0}^n \binom{n}{j} C(K_{m-1,n-j}).
\]
\item[3)]
$a_1$ occurs with an $i$-set $A_0$ of $A-\{a_1\}$, for some 
$i \ge 1$.  Then there must also be a nonempty subset $B_0$ of
$B$ included, say a $j$-set of $B$ with $j \ge 1$.  After 
$A_0$ and $B_0$ are chosen, the remaining elements can be
associated in $C(K_{m-1-i,n-j})$ ways.  The total in this case
is 
\[
\sum_{i=1}^{m-1} \binom{m-1}{i} \sum_{j=1}^n \binom{n}{j}
C(K_{m-1-i,n-j}).
\]
\end{trivlist}
Putting the cases together, we have
\begin{equation}
C(K_{m,n})=
\sum_{j=0}^n \binom{n}{j} C(K_{m-1,n-j})+
\sum_{i=1}^{m-1} \sum_{j=1}^n \binom{m-1}{i} \binom{n}{j}
C(K_{m-1-i,n-j}). \label{Kmnrecurrence}
\end{equation}
Rewrite this as
\begin{equation}
C(K_{m,n})=
\sum_{i=0}^{m-1} \sum_{j=0}^n \binom{m-1}{i} \binom{n}{j}
C(K_{m-1-i,n-j})-
\sum_{i=0}^{m-1} \binom{m-1}{i} C(K_{m-1-i,n}).\label{Kmnbetter}
\end{equation}
Now we can establish that sums of powers of successive
integers arise by induction.  First note
\begin{equation}
\sum_{j=0}^n  \binom{n}{j} C(K_{m-1-i,n-j})=
\sum_{j=0}^n  \binom{n}{j} \sum_{k=1}^{m-i}a_{m-1-i,k} k^{n-j}=
\sum_{k=1}^{m-i} a_{m-1-i,k} (k+1)^n,
\end{equation}
which repeatedly uses the identity 
\[\sum_{j=0}^n \binom{n}{j} x^n =(x+1)^n.\]
The proof is completed by equating coefficients of $k^n$ in (\ref{Kmnbetter}).  Padding the table of coefficients with an
initial column of 0s makes the recurrence work unaltered for
$a_{m,1}$.
$\Box$

Here is a brief table of the coefficients $a_{i,j}$ that the
binomial coefficient summations produce.\newline  
$\begin{array}{rrrrrrrrrr} 
n\backslash i&
         1&  2&   3&   4&  5&  6& 7&8 &9   \\
0\mid&   1                                 \\
1\mid&   0&  1                             \\
2\mid&  -1&  1&   1                        \\
3\mid&  -1& -2&   3&    1                  \\
4\mid&   2& -9&   1&    6&   1             \\
5\mid&   9& -9& -25&   15&  10&  1         \\ 
6\mid&   9& 50&-104&  -20&  50& 15& 1      \\
7\mid& -50&267& -98& -364& 105&119& 21& 1  \\
8\mid&-267&413&1163&-1610&-539&574&238&28&1
\end{array}$

Several properties
of this array follow from the series expansion:
\begin{trivlist}
\item[1.] The main diagonal entry is always 1.
\item[2.] The second diagonal consists of triangular numbers.
\item[3.] Further diagonals are values of polynomials in $n$ as well.
The next three diagonals are represented by polynomials of degrees
4, 6, and 8.
\item[4.] The row sum of each row is 1.
\item[5.] The alternating row sum of each row, taking the main
diagonal entry as positive, is 1.
\item[6.] The first two columns have values that match, up to
a shift and change of sign.  The first column consists of 
coefficients of the series expansion of $e^{1-e^x}$.
\end{trivlist}

This last property is perhaps more than moonshine, given the
generating function of $B(n)$ and the inclusion of all edges 
(subject to one constraint) in $K_{m,n}$.

A few values of $C(K_{m,n})$ calculated from (\ref{a}) and 
(\ref{Ksum}) are given below.\newline
\small{$\begin{array}{rrrrrrrrr}\\
m\backslash n&
      1&2&3&4&5&6&7&8\\
1 \mid & 2& 4& 8& 16& 32& 64& 128& 256\\ 
2 \mid & 4& 12& 34& 96& 274& 792& 2314& 6816\\
3 \mid & 8& 34& 128& 466& 1688& 6154& 22688& 84706\\
4 \mid & 16& 96& 466& 2100& 9226& 40356& 177466& 788100\\ 
5 \mid & 32& 274& 1688& 9226& 48032& 245554& 1251128& 6402586\\ 
6 \mid & 64& 792& 6154& 40356& 245554& 1444212& 8380114& 48510036\\
7 \mid & 128& 2314& 22688& 177466& 1251128& 8380114& 54763088& 
354298186
\end{array}$}
 
%Cubes $Q_n$

%Ladders that wrap around themselves like cycles

\section{Prospectus}

There are several directions that we expect further work on
graph compositions to take.  First, there are many other families
of graphs that have been studied in the literature, and at least 
some of them seem to be appropriate to analyse in the manner of 
this paper.  

The algorithms we have developed to count (and
represent in diagrams) graph compositions are 
sufficiently efficient to handle graphs with up to 20 edges, so 
that, for
instance, we can calculate that the Petersen graph has exactly
8581 compositions. This is important for this paper, if for
no other reason because
every paper in graph theory should mention the Petersen
graph at least once.  Extended numerical data awaits the development of more efficient algorithms.

Another project is to develop a calculus of graph compositions,
so that, for example, we can predict how the number of
compositions is affected when two disjoint graphs
are joined by $k$ edges, or when one or more (adjacent or
nonadjacent) edges are deleted from a given
graph.  Theorems \ref{disconnected}, \ref{one edge}, and 
\ref{complete-1} are small steps in this direction.  We would
like to say something about how operations such as union, 
product, or join of graphs combine the number of compositions.  
\cite{kmr} develops some more tools and uses them to analyze another class of graphs.

% endmatter
\vspace{.1in}
{\small
\begin{flushleft}
Department of Computational and Applied Mathematics\\
Wits 2050, Johannesburg, South Africa\\
{\tt arnoldk@gauss.cam.wits.ac.za}\\
{\tt http://sunsite.wits.ac.za/wits/science/number\_theory/apublic.htm}\\
\vspace{.1in}
Department of Mathematics\\ 
West Virginia University, Morgantown WV 26506-6310\\
{\tt mays@math.wvu.edu}\\
{\tt http://www.math.wvu.edu/homepages/mays}
\end{flushleft}
}

\end{document}